\documentclass[12pt]{amsart}
\newtheorem{theorem}{Theorem}
\newtheorem{corollary}{Corollary}

\newtheorem{lemma}{Lemma}
\begin{document}
\title{On the curvature groups of a CR manifold}
\author{Elisabetta BARLETTA \and Sorin DRAGOMIR}
\maketitle

\centerline{\small\em Dedicated to the memory of Aldo Cossu}

\begin{abstract} We show that any contact form whose Fefferman metric admits a nonzero
parallel vector field is pseudo-Einstein of constant
pseudohermitian scalar curvature. As an application we compute the
curvature groups $H^k (C(M), \Gamma )$ of the Fefferman space
$C(M)$ of a strictly pseudoconvex real hypersurface $M \subset
{\mathbb C}^{n+1}$.
\end{abstract}

\section{Statement of results}
Let $M$ be a strictly pseudoconvex CR manifold of CR dimension $n$
and $\theta$ a contact form on $M$ such that the Levi form
$L_\theta$ is positive definite. Let $S^1 \to C(M) \to M$ be the
canonical circle bundle and $F_\theta$ the Fefferman metric on
$C(M)$, cf. \cite{kn:Lee1}. Let ${\rm GL}(2n+2, {\mathbb R}) \to
L(C(M)) \to C(M)$ be the principal bundle of linear frames tangent
to $C(M)$ and $\Gamma : u \in L(C(M)) \mapsto \Gamma_u \subset T_u
(L(C(M)))$ the Levi-Civita connection of $F_\theta$. Let $H^k
(C(M), \Gamma )$ be the curvature groups
 of $(C(M), \Gamma )$, cf. \cite{kn:GoPe} and
our Section 2. Our main result is
\begin{theorem} If $(M , \theta )$ is a pseudo-Einstein manifold
of constant pseudohermitian scalar curvature $\rho$ then the
curvature groups $H^k (C(M), \Gamma )$ are isomorphic to the de
Rham cohomology groups of $C(M)$. Otherwise {\rm (}that is if
either $\theta$ is not pseudo-Einstein or $\rho$ is
nonconstant{\rm )} $H^k (C(M), \Gamma ) = 0$, $1 \leq k \leq
2n+2$. \label{t:1}
\end{theorem}
The key ingredient in the proof of Theorem \ref{t:1} is the
explicit calculation of the infinitesimal conformal
transformations of the Lorentz manifold $(C(M), F_\theta )$.
\begin{corollary} Let $\Omega \subset {\mathbb C}^{n+1}$ be a smoothly
bounded strictly pseudoconvex domain. There is a defining function
$\varphi$ of $\Omega$ such that $\theta =
\frac{i}{2}(\overline{\partial} -
\partial ) \varphi$ is a pseudo-Einstein contact form on $\partial \Omega$.
If $(\partial \Omega , \theta )$ has constant
pseudohermitian scalar curvature then
\[ H^k (C(\partial \Omega ) , \Gamma ) \approx H^k (\partial \Omega , {\mathbb R}) \oplus
H^{k-1}(\partial \Omega , {\mathbb R}), \] for any $1 \leq k \leq
2n+2$. \label{c:1}
\end{corollary}
The first statement in Corollary \ref{c:1} is a well known
consequence of the fact that $T_{1,0}(\partial \Omega )$ is an
embedded CR structure, cf. J.M. Lee, \cite{kn:Lee2}. If for
instance $\Omega$ is the unit ball in ${\mathbb C}^{n+1}$ and
$\theta = \frac{i}{2}(\overline{\partial} - \partial )|z|^2$ then
\[ H^k (C(\partial \Omega ), \Gamma ) = \begin{cases} {\mathbb R}, & k \in
\{ 1, 2n+1, 2n+2 \} , \cr 0, & {\rm otherwise}. \cr \end{cases} \]
The paper is organized as follows. In Section 2 we recall S.I.
Goldberg \& N.C. Petridis' curvature groups of a torsion-free
linear connection (cf. also I. Vaisman, \cite{kn:Vai}) as well as
the needed material on CR manifolds, Tanaka-Webster connection and
the Fefferman metric. Section 3 is devoted to the proof of Theorem
\ref{t:1} and corollaries.

\section{The curvature groups of the Fefferman metric}
Let $(M, T_{1,0}(M))$ be a $(2n+1)$-dimensional connected strictly
pseudoconvex CR manifold with the CR structure $T_{1,0}(M) \subset
T(M) \otimes {\mathbb C}$. Let $\theta$ be a contact form on $M$
such that the Levi form
\[ L_\theta (Z, \overline{W}) = - i (d \theta )(Z , \overline{W}),
\;\;\; Z, W \in T_{1,0}(M), \] is positive definite. Let $H(M) =
{\rm Re} \{ T_{1,0}(M) \oplus T_{0,1}(M) \}$ be the Levi
distribution and
\[ J : H(M) \to H(M), \;\;\; J(Z + \overline{Z}) = i(Z -
\overline{Z}), \;\;\; Z \in T_{1,0}(M), \] its complex structure.
Let ${\mathbb C} \to K(M) \to M$ be the complex line bundle
\[ K(M)_x = \{ \omega \in \Lambda^{n+1} T^*_x (M) \otimes {\mathbb
C} : T_{0,1}(M)_x \, \rfloor \, \omega = 0 \} , \;\;\; x \in M. \]
There is a natural action of ${\mathbb R}_+$ (the multiplicative
positive reals) on $K(M) \setminus \{ 0 \}$ such that $C(M) =
(K(M) \setminus \{ 0 \})/{\mathbb R}_+$ is a principal
$S^1$-bundle $\pi : C(M) \to M$ (the canonical circle bundle over
$M$, cf. \cite{kn:DrTo}, Chapter 2). The Fefferman metric
$F_\theta$ is given by
\begin{equation}
F_\theta = \pi^* \tilde{G}_\theta + 2 (\pi^* \theta ) \odot \sigma
, \label{e:1}
\end{equation}
\begin{equation}
\sigma = \frac{1}{n+2}\left\{ d \gamma + \pi^* \left( i\,
{\omega_\alpha}^\alpha -\frac{i}{2}\; g^{\alpha\overline{\beta}} d
g_{\alpha\overline{\beta}} -\frac{\rho}{4(n+1)} \;
\theta\right)\right\} . \label{e:2}
\end{equation}
The Fefferman metric is a Lorentz metric on $C(M)$, cf. J.M. Lee,
\cite{kn:Lee1}. The following conventions are adopted as to the
formulae (\ref{e:1})-(\ref{e:2}). Let $T$ be the characteristic
direction of $d \theta$ i.e. the tangent vector field on $M$
determined by $\theta (T) = 1$ and $T \, \rfloor \, d \theta = 0$.
We set
\[ G_\theta (X,Y) = (d \theta )(X, J Y), \;\;\; X,Y \in H(M), \]
\[ \tilde{G}_\theta (X,Y) = G_{\theta}(X, Y), \;\;\; \tilde{G}_\theta (Z, T) = 0, \;\;\; Z \in T(M). \]
There is a unique linear connection $\nabla$ on $M$
(the Tanaka-Webster connection of $(M , \theta )$, cf.
\cite{kn:Tan} and \cite{kn:Web}) such that i) the Levi
distribution is parallel with respect to $\nabla$, ii) $\nabla J =
0$ and $\nabla g_\theta = 0$, iii) the torsion $T_\nabla$ of
$\nabla$ is pure i.e.
\[ T_\nabla (Z,W) =0, \;\;\; T_\nabla (Z , \overline{W}) = 2 i
L_\theta (Z, \overline{W})T, \;\;\; Z, W \in T_{1,0}(M), \]
\[ \tau \circ J + J \circ \tau = 0. \]
Here $g_\theta$ is the Webster metric i.e. the Riemannian metric
on $M$ given by
\[ g_\theta (X,Y) = G_\theta (X,Y), \;\;\; g_\theta (X,T) = 0,
\;\;\; g_\theta (T,T) = 1, \] for any $X,Y \in H(M)$. Also $\tau
(X) = T_\nabla (T, X)$, $X \in T(M)$, is the pseudohermitian
torsion of $\nabla$. If $\{ T_\alpha : 1 \leq \alpha \leq n \}$ is
a local frame of $T_{1,0}(M)$ defined on the open set $U \subseteq
M$ then ${\omega_\alpha}^\beta$ are the corresponding connection
$1$-forms of the Tanaka-Webster connection i.e. $\nabla T_\alpha =
{\omega_\alpha}^\beta \otimes T_\beta$. Let $R^\nabla$ be the
curvature of $\nabla$ and
\[ R_{\alpha\overline{\beta}} = {\rm trace} \{ X \mapsto R^\nabla
(X, T_\alpha ) T_{\overline{\beta}} \} \] the pseudohermitian
Ricci tensor of $(M , \theta )$. Moreover
$g_{\alpha\overline{\beta}} = L_\theta (T_\alpha ,
T_{\overline{\beta}})$ and $\rho = g^{\alpha\overline{\beta}}
R_{\alpha\overline{\beta}}$ is the pseudohermitian scalar
curvature of $\nabla$. Also $\gamma : \pi^{-1}(U) \to {\mathbb R}$
is a local fibre coordinate on $C(M)$. Precisely let $\{
\theta^\alpha : 1 \leq \alpha \leq n \}$ be the admissible local
coframe associated to $\{ T_\alpha : 1 \leq \alpha \leq n \}$ i.e.
\[ \theta^\alpha (T_\beta ) = \delta^\alpha_\beta \, , \;\;\;
\theta^\alpha (T_{\overline{\beta}}) = 0, \;\;\; \theta^\alpha (T)
= 0. \] The locally trivial structure of $S^1 \to C(M) \to M$ is
described by
\[ \pi^{-1}(U) \to U \times S^1 \, , \;\;\; [\omega ] \mapsto (x
\, , \, \frac{\lambda}{|\lambda |}), \;\;\; \omega \in K(M)_x
\setminus \{ 0 \} , \]
\[ \omega = \lambda (\theta \wedge \theta^1 \wedge \cdots \wedge
\theta^n )_x \, , \;\;\; x \in U, \;\;\; \lambda \in {\mathbb C}
\setminus \{ 0 \} . \] Then $\gamma ([\omega ]) = \arg (\lambda
/|\lambda |)$ where $\arg : S^1 \to [0, 2 \pi )$. If $(U, x^1 ,
\cdots , x^{2n+1})$ is a system of local coordinates on $M$ then
$(\pi^{-1}(U), \tilde{x}^1 , \cdots , \tilde{x}^m )$ are the
naturally induced local coordinates on $C(M)$ i.e. $\tilde{x}^A =
x^A \circ \pi$, $1 \leq A \leq 2n+1$, and $\tilde{x}^m = \gamma$
(with $m = 2n+2$).
\par
Let $\Pi : L(C(M)) \to C(M)$ be the projection and $\rho : {\rm
GL}(m, {\mathbb R}) \to {\rm End}_{\mathbb R}({\mathbb R}^m )$ the
natural representation. We denote by $\Omega^k_{\rho ({\rm
GL}(m))}(C(M))$ the space of tensorial $k$-forms of type $\rho
({\rm GL}(m, {\mathbb R}))$ i.e. each $\omega \in \Omega^k_{\rho
({\rm GL}(m))} (C(M))$ is a ${\mathbb R}^m$-valued $k$-form on
$L(C(M))$ such that
\par
i) $\omega_u (X_1 , \cdots , X_k ) = 0$ if at least one $X_i \in
{\rm Ker}(d_u \Pi )$,
\par
ii) $\omega_{u g} ((d_u R_g ) X_1 , \cdots , (d_u R_g ) X_k ) =
\rho (g^{-1}) \omega_u (X_1 , \cdots , X_k )$ for any $g \in {\rm
GL}(m , {\mathbb R})$, $X_i \in T_u (L(C(M)))$ and $u \in
L(C(M))$.
\par
Let $\Gamma$ be the Levi-Civita connection of $(C(M), F_\theta )$
thought of as a connection-distribution in $L(C(M)) \to C(M)$. If
$\omega \in \Omega^k_{\rho ({\rm GL}(m))}(C(M))$ then its
covariant derivative with respect to $\Gamma$ is the tensorial
$(k+1)$-form of type $\rho ({\rm GL}(m, {\mathbb R}))$
\[ (\nabla \omega )(X_0 , \cdots , X_k ) = \omega (h X_0 , \cdots
, h X_k ), \] for any $X_i \in T(L(C(M)))$, $0 \leq i \leq k$.
Here $h_u : T_u (L(C(M))) \to \Gamma_u$ is the natural projection
associated to the direct sum decomposition $T_u (L(C(M))) =
\Gamma_u \oplus {\rm Ker}(d_u \Pi )$. Let us consider the
$C^\infty (C(M))$-module
\[ L^k = \Omega^k_{\rho ({\rm GL}(m))} (C(M)) \times \Pi^*
\Omega^k (C(M)) \] and the submodule $\tilde{L}^k$ given by
\[ \tilde{L}^k = \{ (\omega , \Pi^* \alpha ) \in L^k : \nabla^2
\omega = 0 \} . \] Let $\eta \in \Gamma^\infty (T^* (L(C(M)))
\otimes {\mathbb R}^m )$ be the canonical $1$-form i.e. $\eta_u =
u^{-1} \circ (d_u \Pi )$ for any $u \in L(C(M))$. If we set
\[ D^k : \tilde{L}^k \to \tilde{L}^{k+1} , \;\;\; D^k (\omega , \Pi^* \alpha ) = (\nabla \omega - \eta \wedge
\Pi^* \alpha \, , \, \Pi^* d \alpha ), \] then $\tilde{L} =
(\bigoplus_{k=0}^m \tilde{L}^k \, , \, D^k )$ is a cochain
complex, cf. \cite{kn:GoPe}, p. 550. The curvature groups of
$\Gamma$ are the cohomology groups
\[ H^k (C(M), \Gamma ) = H^k (\tilde{L}) = \frac{{\rm Ker}(D^k
)}{D^{k-1} \, \tilde{L}^{k-1}} \, , \;\;\; 1 \leq k \leq m. \]

\section{Infinitesimal conformal transformations}
We shall establish the following
\begin{theorem} Any infinitesimal conformal transformation of
the Fefferman metric $F_\theta$ is a parallel vector field.
\label{t:2}
\end{theorem}
By a result of C.R. Graham, \cite{kn:Gra}, $\sigma$ is a
connection $1$-form in $S^1 \to C(M) \to M$. For each vector field
$X \in T(M)$ let $X^\uparrow$ denote the horizontal lift of $X$
with respect to $\sigma$ i.e. $X^\uparrow_z \in {\rm Ker}(\sigma_z
)$ and $(d_z \pi ) X^\uparrow_z = X_{\pi (z)}$ for any $z \in
C(M)$. To prove Theorem \ref{t:2} we need to recall the following
\begin{lemma} {\rm (E. Barletta et al., \cite{kn:BaDrUr})} \par\noindent
The Levi-Civita connection $\nabla^{C(M)}$ of $(C(M), F_\theta )$
and the Tanaka-Web\-ster connection $\nabla$ of $(M , \theta )$
are related by
\begin{equation}
\nabla^{C(M)}_{X^\uparrow} Y^\uparrow = (\nabla_X Y)^\uparrow  -
(d \theta )(X,Y) T^\uparrow - \{ A(X,Y) + (d \sigma )(X^\uparrow ,
Y^\uparrow ) \} S , \label{e:3}
\end{equation}
\begin{equation}
\nabla^{C(M)}_{X^\uparrow} T^\uparrow = (\tau (X) + \phi
X)^\uparrow \, , \label{e:4}
\end{equation}
\begin{equation}
\nabla^{C(M)}_{T^\uparrow} X^\uparrow = (\nabla_T X + \phi
X)^\uparrow + 2 (d \sigma )(X^\uparrow , T^\uparrow ) S,
\label{e:5}
\end{equation}
\begin{equation}
\nabla^{C(M)}_{X^\uparrow} S = \nabla^{C(M)}_S X^\uparrow = (J
X)^\uparrow \, , \label{e:6}
\end{equation}
\begin{equation}
\nabla^{C(M)}_{T^\uparrow} T^\uparrow = V^\uparrow \, , \;\;\;
\nabla^{C(M)}_S S = 0, \label{e:7}
\end{equation}
\begin{equation}
\nabla^{C(M)}_S T^\uparrow = 0, \;\;\; \nabla^{C(M)}_{T^\uparrow}
S = 0, \label{e:8}
\end{equation}
for any $X,Y \in H(M)$. Here $A(X,Y) = g_\theta (\tau (X), Y)$.
Also the vector field $V \in H(M)$ and the endomorphism $\phi :
H(M) \to H(M)$ are given by
\[ G_\theta (V , Y) = 2 (d \sigma )(T^\uparrow , Y^\uparrow ), \;\;\;
G_\theta (\phi X \, , \, Y) = 2 (d \sigma )(X^\uparrow ,
Y^\uparrow ), \] for any $X,Y \in H(M)$. \label{l:1}
\end{lemma}
A vector field $\mathcal X$ on $C(M)$ is an infinitesimal
conformal transformation of $F_\theta$ if
\begin{equation}
\nabla^{C(M)} {\mathcal X} = \lambda I, \label{e:9}
\end{equation}
for some $\lambda \in C^\infty (C(M))$ where $I$ is the identical
transformation of $T(C(M))$. Let $S$ be the tangent to the
$S^1$-action (locally $S = \partial /\partial \gamma$). Taking
into account the decomposition $T(C(M)) = H(M)^\uparrow \oplus
{\mathbb R} T^\uparrow \oplus {\mathbb R} S$ the first order
partial differential system (\ref{e:9}) is equivalent to
\begin{equation}
\nabla^{C(M)}_{X^\uparrow} {\mathcal X} = \lambda X^\uparrow \, ,
\;\; \nabla^{C(M)}_{T^\uparrow} {\mathcal X} = \lambda T^\uparrow
\, , \;\; \nabla^{C(M)}_S {\mathcal X} = \lambda S,
\label{e:10}
\end{equation}
for any $X \in H(M)$. Let $\{ X_a : 1 \leq a \leq 2n \} = \{
X_\alpha , \; J X_\alpha : 1 \leq \alpha \leq n \}$ be a local
frame of $H(M)$. Then ${\mathcal X} = {\mathcal X}^a X_a^\uparrow
+ f T^\uparrow + g S$ for some $C^\infty$ functions ${\mathcal
X}^a$, $f$ and $g$. By (\ref{e:6})-(\ref{e:8}) in Lemma \ref{l:1}
the last equation in (\ref{e:10}) may be written
\[ S({\mathcal X}^a ) X_a^\uparrow + {\mathcal X}^a
(J X_a )^\uparrow + S(f) T^\uparrow + S(g) S = \lambda S
\]
hence
\begin{equation} {\mathcal X}^a = 0, \;\;\; S(f) = 0, \;\;\; S(g)
= \lambda . \label{e:11}
\end{equation}
Similarly, by (\ref{e:4}) and (\ref{e:6}) in Lemma \ref{l:1} the
first equation in (\ref{e:10}) may be written
\[ X^\uparrow (f) T^\uparrow + f(\tau (X) + \phi X)^\uparrow +
X^\uparrow (g) S + g (J X )^\uparrow = \lambda X^\uparrow \] hence
\begin{equation} X^\uparrow (f) = 0, \;\;\; X^\uparrow (g) = 0,
\label{e:12}
\end{equation}
\begin{equation}
f(\tau (X) + \phi X)^\uparrow + g (J X)^\uparrow = \lambda
X^\uparrow . \label{e:13}
\end{equation}
\begin{lemma} With respect to a local frame $\{ T_\alpha : 1 \leq \alpha \leq n
\}$ of $T_{1,0}(M)$ the endomorphism $\phi : H(M) \otimes {\mathbb
C} \to H(M) \otimes {\mathbb C}$ is given by $\phi T_\alpha =
{\phi_\alpha}^\beta T_\beta + {\phi_\alpha}^{\overline{\beta}}
T_{\overline{\beta}}$ with
\begin{equation}
\phi^{\overline{\alpha}\beta} = \frac{i}{2(n+2)} \{
R^{\overline{\alpha}\beta} - \frac{\rho}{2(n+1)} \;
g^{\overline{\alpha}\beta} \} , \;\;\; \phi^{\alpha\beta} = 0,
\label{e:14}
\end{equation}
and $\phi^{\overline{\alpha}\beta} = g^{\overline{\alpha}\gamma}
{\phi_\gamma}^\beta$ and $\phi^{\alpha\beta} =
g^{\alpha\overline{\gamma}} = {\phi_{\overline{\gamma}}}^\beta$.
\label{l:2}
\end{lemma}
{\em Proof of Lemma} \ref{l:2}. Taking the exterior derivative of
(\ref{e:2}) we obtain
\[ (n+2) d \sigma = \pi^* \left( i d {\omega_\alpha}^\alpha -
\frac{i}{2} d g^{\alpha\overline{\beta}} \wedge d
g_{\alpha\overline{\beta}} - \frac{1}{4(n+1)} \, d (\rho \theta )
\right) . \] Note that $\nabla g_\theta = 0$ may be locally
written as $d g_{\alpha\overline{\beta}} =
g_{\alpha\overline{\gamma}}
{\omega_{\overline{\beta}}}^{\overline{\gamma}} +
{\omega_\alpha}^\gamma g_{\gamma\overline{\beta}}$. Also
$g^{\alpha\overline{\beta}} g_{\overline{\beta}\gamma} =
\delta^\alpha_\gamma$ yields $d g^{\alpha\overline{\beta}} = -
g^{\gamma\overline{\beta}} g^{\alpha\overline{\rho}} d
g_{\overline{\rho}\gamma}$. Hence
\[ d g^{\alpha\overline{\beta}} \wedge d
g_{\alpha\overline{\beta}} = \omega_{\alpha\overline{\beta}}
\wedge \omega^{\alpha\overline{\beta}} +
\omega_{\overline{\alpha}\beta} \wedge
\omega^{\overline{\alpha}\beta} = 0. \] Let $\{ \theta^\alpha : 1
\leq \alpha \leq n \}$ be the admissible local coframe associated
to $\{ T_\alpha : 1 \leq \alpha \leq n \}$. Then (by a result in
\cite{kn:Web}, cf. also \cite{kn:DrTo}, Chapter 1)
\[ d {\omega_\alpha}^\alpha = R_{\lambda\overline{\mu}}
\theta^\lambda \wedge \theta^{\overline{\mu}} +
(W^\alpha_{\alpha\lambda} \theta^\lambda -
W^\alpha_{\alpha\overline{\mu}} \theta^{\overline{\mu}}) \wedge
\theta , \]
\[ W^\beta_{\alpha\lambda} = g^{\overline{\sigma}\beta} \nabla_{\overline{\sigma}} A_{\alpha\lambda} \, ,
\;\;\; W_{\alpha\overline{\mu}}^\beta = g^{\overline{\sigma}\beta}
\nabla_\alpha A_{\overline{\mu}\; \overline{\sigma}} \, , \] where
$A_{\alpha\beta} = A(T_\alpha , T_\beta )$ and covariant
derivatives are meant with respect to the Tanaka-Webster
connection. Finally (by the very definition of $\phi$)
\[ (n+2) G_\theta (\phi X , Y) = i (R_{\alpha\overline{\beta}}
\theta^\alpha \wedge \theta^{\overline{\beta}} )(X,Y) -
\frac{\rho}{4(n+1)} (d \theta )(X,Y) \] for any $X,Y \in H(M)
\otimes {\mathbb C}$. This yields (\ref{e:14}). Lemma \ref{l:2} is
proved.
\par
{\em Proof of Theorem} \ref{t:2}. Let us extend both members of
(\ref{e:13}) by ${\mathbb C}$-linearity. Then (\ref{e:13}) holds
for any $X \in H(M) \otimes {\mathbb C}$. By a result in
\cite{kn:Tan} $\tau (T_{1,0}(M)) \subseteq T_{0,1}(M)$ hence $\tau
(T_\alpha ) = A_\alpha^{\overline{\beta}} T_{\overline{\beta}}$
for some $C^\infty$ functions $A_\alpha^{\overline{\beta}}$. Using
(\ref{e:13}) for $X = T_\alpha$ we obtain
\begin{equation}
f \, A_\alpha^{\overline{\beta}} = 0, \;\;\; f \,
{\phi_\alpha}^\beta + (i g - \lambda ) \delta_\alpha^\beta = 0.
\label{e:15}
\end{equation}
By Lemma \ref{l:2}
\[ {\phi_\alpha}^\beta = \frac{i}{2(n+2)} \left( {R_\alpha}^\beta
- \frac{\rho}{2(n+1)} \; \delta_\alpha^\beta \right) \] and a
contraction leads to ${\phi_\alpha}^\alpha = i \rho /[4(n+1)]$.
Next a contraction in the second of the identities (\ref{e:15})
gives $f \; {\phi_\alpha}^\alpha + n(i g - \lambda ) = 0$ or $i
\rho f + 4n(n+1) i g = 4n(n+1) \lambda$ and then
\begin{equation}
g = - \frac{\rho}{4n(n+1)} \; f \label{e:16}
\end{equation}
and $\lambda = 0$ as $f$, $g$ and $\lambda$ are ${\mathbb
R}$-valued. In particular $\nabla^{C(M)} {\mathcal X} = 0$.
Theorem \ref{t:2} is proved.
\begin{corollary} Any infinitesimal conformal transformation of
$F_\theta$ is a parallel vector field of the form
\begin{equation}
{\mathcal X} = a \left( T^\uparrow - \frac{\rho}{4n(n+1)} \; S
\right) , \;\;\; a \in {\mathbb R}. \label{e:17}
\end{equation}
In particular any contact form $\theta$ whose Fefferman metric
$F_\theta$ admits a nontrivial parallel vector field is
pseudo-Einstein of constant pseudohermitian scalar curvature and
vanishing pseudohermitian torsion. \label{c:2}
\end{corollary}
{\em Proof}. By (\ref{e:7})-(\ref{e:8}) in Lemma \ref{l:1} the
middle equation in (\ref{e:10}) may be written
\[ T^\uparrow (f) T^\uparrow + f V^\uparrow + T^\uparrow (g) S =
\lambda T^\uparrow \] hence
\begin{equation}
T^\uparrow (f) = \lambda , \;\;\; T^\uparrow (g) = 0, \label{e:18}
\end{equation}
\begin{equation}
f(z) \; V_{\pi (z)} = 0, \;\;\; z \in C(M). \label{e:19}
\end{equation}
Yet $\lambda = 0$ (by Theorem \ref{t:2}) so that (by
(\ref{e:11})-(\ref{e:12}) and (\ref{e:18})) $f = a$ and $g = b$
for some $a,b \in {\mathbb R}$. Let us assume now that $F_\theta$
admits a parallel vector field ${\mathcal X} \neq 0$. Replacing
from (\ref{e:16}) into the second of the identities (\ref{e:15})
leads to
\begin{equation} a \left( {R_\alpha}^\beta - \frac{\rho}{n} \,
\delta_\alpha^\beta \right) = 0. \label{e:20}
\end{equation}
Note that $a \neq 0$ (otherwise (\ref{e:16}) implies $b = 0$ hence
${\mathcal X} = 0$) so that (by (\ref{e:20}))
$R_{\alpha\overline{\beta}} = (\rho /n)
g_{\alpha\overline{\beta}}$ i.e. $\theta$ is pseudo-Einstein (cf.
\cite{kn:Lee2}). Also (\ref{e:16}) shows that $\rho =$ constant.
Finally the first identity in (\ref{e:15}) implies $\tau = 0$.
Corollary \ref{c:2} is proved.
\par
A remark is in order. Apparently (\ref{e:19}) implies that $a = 0$
when ${\rm Sing}(V) \neq \emptyset$ (and then there would be no
nonzero parallel vector fields on $(C(M), F_\theta )$). Yet we may
show that
\begin{corollary} Assume that $\theta$ is pseudo-Einstein. Then $V = 0$ if and only if $\rho$ is
constant. \label{c:3}
\end{corollary}
So (\ref{e:19}) brings no further restriction. {\em Proof of
Corollary} \ref{c:3}. Note that
\[ 2 (d{\omega_\alpha}^\alpha )(T, T_\beta ) = -
W^\alpha_{\alpha\beta} \, , \;\;\; 2 (d{\omega_\alpha}^\alpha )(T,
T_{\overline{\beta}}) = W^\alpha_{\alpha\overline{\beta}} \, , \]
\[ 2 d(\rho \theta )(T, T_\beta ) = - \rho_\beta \, , \;\;\; 2 d
(\rho \theta )(T, T_{\overline{\beta}}) = -
\rho_{\overline{\beta}} \, , \] where $\rho_\beta = T_\beta (\rho
)$ and $\rho_{\overline{\beta}} = \overline{\rho_\beta}$.
Consequently
\[ 2(n+2) (d \sigma )(T^\uparrow , T_\beta^\uparrow ) = - i
W^\alpha_{\alpha\overline{\beta}} + \frac{1}{4(n+1)} \, \rho_\beta
\, .
\] On the other hand (cf. \cite{kn:DrTo}, Chapter 5) if $\theta$
is pseudo-Einstein then
\[ W^\alpha_{\alpha\beta} = - \frac{i}{2n} \; \rho_\beta \, ,
\;\;\; W^\alpha_{\alpha\overline{\beta}} =
\overline{W^\alpha_{\alpha\beta}} \, . \] Hence $V$ is given by
(see Lemma \ref{l:1} above)
\[ G_\theta (V , T_\beta ) = - \frac{1}{4n(n+1)} \; \rho_\beta \, .  \]
Clearly if $\rho =$ constant then $V = 0$. Conversely if $V = 0$
then $\overline{\partial}_b \rho = 0$ i.e. $\rho$ is a ${\mathbb
R}$-valued CR function. As $M$ is nondegenerate $\rho$ is
constant. Corollary \ref{c:3} is proved.
\par
At this point we may prove Theorem \ref{t:1}. Let ${\mathcal L}^k$
be the sheaf associated to the module $\tilde{L}^k$  i.e. for any
open set $A \subseteq C(M)$
\[ {\mathcal L}^k (A) = \{ (\lambda , \Pi^* \alpha ) : \lambda \in
\Omega^k_{\rho ({\rm GL}(m))}(\Pi^{-1}(A)), \;\; \nabla^2 \lambda
= 0, \;\; \alpha \in \Omega^k (A)\} . \] Let $D^k : {\mathcal L}^k
\to {\mathcal L}^{k+1}$ be the sheaf homomorphism induced by the
module homomorphism $D^k : \tilde{L}^k \to \tilde{L}^{k+1}$.
\begin{lemma} Let ${\mathcal S}_\theta$ be the sheaf of parallel vector fields on
$(C(M), F_\theta )$. For each open set $A \subseteq C(M)$ let $j_A
: {\mathcal S}_\theta (A) \to {\mathcal L}^0 (A)$be given by \[
j_A ({\mathcal X}) = (f_{\mathcal X} , 0), \;\;\; {\mathcal X} \in
{\mathcal S}_\theta (A),
\]
\[ f_{\mathcal X} : \Pi^{-1}(A) \to {\mathbb R}^m \, , \;\;\;
f_{\mathcal X}(u) = u^{-1} ({\mathcal X}_{\Pi (u)} ), \;\;\; u \in
\Pi^{-1} (A). \] Then
\begin{equation}
0 \to {\mathcal S}_\theta \stackrel{j}{\longrightarrow} {\mathcal
L}^0 \stackrel{D^0}{\longrightarrow} {\mathcal L}^1
\stackrel{D^1}{\longrightarrow} \cdots
\stackrel{D^{m-1}}{\longrightarrow} {\mathcal L}^m \to 0
\end{equation}
is a fine resolution of ${\mathcal S}_\theta$ so that the
curvature groups of $\Gamma$ are isomorphic to the cohomology
groups of $C(M)$ with coefficients in ${\mathcal S}_\theta$.
\label{l:3}
\end{lemma}
{\em Proof}. Let $(f , \lambda ) \in {\mathcal L}^0$ such that $0
= D^0 (f, \lambda ) = (\nabla f - \lambda \eta \, , \, d \lambda
)$. Let $z \in C(M)$ and $u \in \Pi^{-1}(z)$. We set by definition
${\mathcal X}_z = u(f(u))$. As $f \circ R_g = \rho (g^{-1}) \circ
f$ for any $g \in {\rm GL}(m , {\mathbb R})$ it follows that
${\mathcal X}_z$ is well defined. Let $(\Pi^{-1}(C(U)), x^i ,
X^i_j )$ be the naturally induced local coordinates on $L(C(M))$
i.e. $x^i (u) = \tilde{x}^i (\Pi (u))$ and $X^i_j (u) = a^i_j$ for
any $u = (z, \{ X_i : 1 \leq i \leq m \} ) \in L(C(M))$ such that
$X_j = a_j^i (\partial /\partial \tilde{x}^i )_z$ (here $C(U) =
\pi^{-1}(U)$). Given a vector field $X = \lambda^j \partial
/\partial x^j + \lambda^i_j \partial /\partial X^i_j$ on $L(C(M))$
\[ (\nabla f)_u X_u = \lambda^j (u) \; u^{-1} (\nabla^{C(M)}_{\partial /\partial \tilde{x}^j}
{\mathcal X} )_{\Pi (u)} \, , \;\;\; u \in \Pi^{-1}(C(U)). \] Then
$\nabla f = \lambda \eta$ implies that $\nabla^{C(M)} {\mathcal X}
= \lambda I$ hence (by Theorem \ref{t:2}) $\lambda = 0$ i.e. $(f,
\lambda ) = j ({\mathcal X})$. Therefore the corresponding
sequence of stalks $0 \to {\mathcal S}_{\theta , z} \to {\mathcal
L}^0_z \to {\mathcal L}^1_z \to \cdots \to {\mathcal L}^m_z \to 0$
is exact at ${\mathcal L}^0_z$ while the exactness at the
remaining terms follows from the Poincar\'e lemma for $D$ as in
\cite{kn:GoPe}, p. 552. Lemma \ref{l:3} is proved. In particular
if $\theta$ is a pseudo-Einstein contact form of constant
pseudohermitian scalar curvature then (by Corollary \ref{c:2})
${\mathcal S}_\theta = {\mathbb R}$ and Lemma \ref{l:3} furnishes
a resolution $0 \to {\mathbb R} \stackrel{j}{\to} {\mathcal L}^*$
of the constant sheaf ${\mathbb R}$ where
\[ j_A (a) = (f_{\mathcal X} , 0), \;\;\; a \in {\mathbb R}, \]
and ${\mathcal X}$ is given by (\ref{e:17}) in Corollary
\ref{c:2}, hence
\[ H^k (C(M), \Gamma ) \approx H^k (C(M), {\mathbb R}), \;\;\; 1
\leq k \leq m. \] Otherwise (i.e. if $\theta$ is not
pseudo-Einstein or $\rho$ is nonconstant) then ${\mathcal
S}_\theta = 0$. Theorem \ref{t:1} is proved.
\par
If $M \subset {\mathbb C}^{n+1}$ is a strictly pseudoconvex real
hypersurface then (by a result of J.M. Lee, \cite{kn:Lee2}) $M$
admits globally defined pseudo-Einstein contact forms. On the
other hand the pullback to $M$ of $d z^0 \wedge \cdots \wedge d
z^n$ is a global nonzero section in $K(M)$. In particular $C(M)$
is trivial. If $\theta$ is a pseudo-Einstein contact form on $M$
of constant pseudohermitian scalar curvature then (by Theorem
\ref{t:1} and the K\"unneth formula)
\[ H^k (C(M), \Gamma ) = H^k (C(M), {\mathbb R} ) = \sum_{p+q=k} H^p (M, {\mathbb R})
\otimes H^q (S^1 , {\mathbb R}) = \]
\[ = H^k (M , {\mathbb R}) \oplus H^{k-1}(M , {\mathbb R}) \]
and Corollary \ref{c:1} is proved. Using M. Rumin's criterion (cf.
\cite{kn:Rum}) for the vanishing of the first Betti number of a
pseudohermitian manifold we get
\begin{corollary} Let $M \subset {\mathbb C}^{n+1}$ be a connected
strictly pseudoconvex real hypersurface and $\theta$ a
pseudo-Einstein contact form on $M$ with $\rho =$ constant and
$\tau = 0$. If $n \geq 2$ then $H^1 (C(M), \Gamma ) = {\mathbb
R}$. \label{c:4}
\end{corollary}
An interesting question is whether one may improve Corollary
\ref{c:1} by choosing a contact form with $\rho$ constant to start
with. Indeed as $T_{1,0}(M)$ is embedded one may choose a
pseudo-Einstein contact form $\theta$. On the other hand if the CR
Yamabe invariant $\lambda (M)$ is $< \lambda (S^{2n+1})$ then (by
the solution to the CR Yamabe problem due to D. Jerison \& J.M.
Lee, \cite{kn:JeLe}) there is a positive solution $u$ to the CR
Yamabe equation such that $u^{2/n} \theta$ has constant
pseudohermitian scalar curvature. Yet (by a result in
\cite{kn:Lee2}) the pseudo-Einstein property is preserved if and
only if $u$ is a CR-pluriharmonic function. It is an open problem
whether the CR Yamabe equation admits CR-pluriharmonic solutions.

\vskip 1cm\noindent {\small Authors' address: Universit\`a degli
Studi della Basilicata, Dipartimento di Matematica, Contrada
Macchia Romana, 85100 Potenza, Italy.
\par\noindent
e-mail: {\tt barletta@unibas.it}, {\tt dragomir@unibas.it}}

\begin{thebibliography}{123456}
\bibitem [1]{kn:BaDrUr} E. Barletta \& S. Dragomir \& H. Urakawa,
{\em Yang-Mills fields on CR manifolds}, submitted to J. Math.
Phys., 2005.

\bibitem [2]{kn:DrTo} S. Dragomir \& G. Tomassini, {\em
Differential geometry and analysis on CR manifolds}, to be
published in Progress in Math., Birkh\"auser, Boston, 2006.

\bibitem [3]{kn:GoPe} S.I. Goldberg \& N.C. Petridis, {\em The
curvature groups of a pseudo-Riemannian manifold}, J. Differential
Geometry, 9(1974), 547-555.

\bibitem [4]{kn:Gra} C.R. Graham, {\em On Sparling's characterization
of Fefferman metrics}, American J. Math., 109(1987), 853-874.

\bibitem [5]{kn:JeLe} D. Jerison \& J.M. Lee, {\em The Yamabe
problem on CR manifolds}, J. Diff. Geometry, 25(1987), 167-197.

\bibitem [6]{kn:Lee1} J.M. Lee, {\em The Fefferman metric and
pseudohermitian invariants}, Trans. A.M.S., (1)296(1986), 411-429.

\bibitem [7]{kn:Lee2} J.M. Lee, {\em Pseudo-Einstein structures on
CR manifolds}, American J. Math., 110(1988), 157-178.

\bibitem [8]{kn:Rum} M. Rumin, {\em Un complexe de formes
diff\'erentielles sur les vari\'et\'es de contact}, C.R. Acad.
Sci. Paris, 330(1990), 401-404.

\bibitem [9]{kn:Tan} N. Tanaka, {\em A differential geometric study
on strongly pseudo-convex manifolds}, Kinokuniya Book Store Co.,
Ltd., Kyoto, 1975.

\bibitem [10]{kn:Vai} I. Vaisman, {\em The curvature groups of a
space form}, Ann. Scuola Norm. Sup. Pisa, 22(1968), 331-341.

\bibitem [11]{kn:Web} S.M. Webster, {\em Pseudohermitian structures
on a real hypersurface}, J. Diff. Geometry, 13(1978), 25-41.

\end{thebibliography}
\end{document}